\documentclass[a4paper]{amsart}

\usepackage[T1]{fontenc}
\usepackage[latin1]{inputenc}
\usepackage{amssymb}
\usepackage{euscript}
\usepackage{amsxtra}
\usepackage{ae}
\usepackage[all]{xy}
\usepackage{enumerate}
\usepackage[lite]{amsrefs}

\theoremstyle{plain}
\newtheorem{thm}{Theorem}
\newtheorem{lemma}[thm]{Lemma}
\newtheorem{prop}[thm]{Proposition}
\newtheorem{cor}[thm]{Corollary}
\theoremstyle{definition}
\newtheorem{defn}[thm]{Definition}
\theoremstyle{remark}
\newtheorem{rmk}[thm]{Remark}
\newtheorem{note}[thm]{Note}

\newcommand{\C}{\mathbb C}
\newcommand{\Z}{\mathbb Z}
\newcommand{\R}{\mathbb R}
\newcommand{\N}{\mathbb N}

\newcommand{\KX}{\mathrm{KX}}
\newcommand{\HX}{\mathrm{HX}}

\newcommand{\vv}{\mathfrak{\bar{c}}}
\newcommand{\vvr}{\mathfrak{\bar{c}^{red}}}
\newcommand{\bv}{\mathfrak{c}}
\newcommand{\bvr}{\mathfrak{c^{red}}}
\newcommand{\Rips}{{\mathcal{P}}}

\newcommand{\dirac}{\mathrm{D}\!\!\!\!/}
\newcommand{\Proj}{{\mathcal{P}}}
\newcommand{\Dirac}{\mathcal{D}}

\newcommand*{\KK}{\mathrm{KK}}
\newcommand*{\RKK}{\mathrm{RKK}}
\newcommand*{\K}{\mathrm{K}}
\newcommand*{\Ktop}{\mathrm{K}^{\mathrm{top}}}

\DeclareMathOperator{\Aut}{Aut}
\DeclareMathOperator{\Res}{Res}
\DeclareMathOperator{\Ind}{Ind}

\DeclareMathOperator{\supp}{supp}

\newcommand*{\Hilm}{{\mathcal{E}}}
\newcommand*{\CI}{{\mathcal{CI}}}
\newcommand*{\CC}{{\mathcal{CC}}}

\newcommand*{\EG}{{\mathcal{E}G}}

\newcommand*{\Mult}{{\mathcal{M}}}
\newcommand*{\Cred}{C^*_{\mathrm{r}}}

\newcommand*{\Comp}{{\mathbb{K}}}
\newcommand*{\Bound}{{\mathbb{B}}}

\newcommand*{\brd}{-\hspace{0pt}}
\newcommand*{\nbd}{\nobreakdash-\hspace{0pt}}

\newcommand*{\abs}[1]{\lvert#1\rvert}
\newcommand*{\ket}[1]{\lvert#1\rangle}
\newcommand*{\bra}[1]{\langle#1\rvert}
\newcommand*{\norm}[1]{\lVert#1\rVert}
\newcommand*{\gen}[1]{\langle#1\rangle}

\newcommand*{\rcross}{\mathbin{\rtimes_{\mathrm{r}}}}
\newcommand*{\cross}{\mathbin{\rtimes}}
\newcommand*{\defeq}{\overset{\mathrm{def}}=}

\begin{document}

\title{A geometric characterization of the Dirac dual Dirac method}
\author{Heath Emerson}
\email{hemerson@math.uni-muenster.de}

\author{Ralf Meyer}
\email{rameyer@math.uni-muenster.de}
\address{Mathematisches Institut\\
         Westf.\ Wilhelms-Universität Münster\\
         Ein\-stein\-stra\-ße 62\\
         48149 Münster\\
         Germany}

\begin{abstract}
  Let~$G$ be a discrete, torsion free group with a finite dimensional
  classifying space~$BG$.  We show that the existence of a
  $\gamma$\nobreakdash-element for such~$G$ is a metric, that is,
  coarse, invariant of~$G$.  We also obtain results for groups with
  torsion.  The method of proof involves showing that a group $G$ 
possesses a $\gamma$-element if and only if a certain coarse (co)-assembly map is an isomorphism.

\end{abstract}


\thanks{This research was supported by the EU-Network \emph{Quantum
  Spaces and Noncommutative Geometry} (Contract HPRN-CT-2002-00280)
  and the \emph{Deutsche Forschungsgemeinschaft} (SFB 478).}

\maketitle

\section{Introduction}
\label{sec:intro}
 
The \emph{Descent Principle} asserts that the strong Novikov
conjecture for a discrete group~$G$ with finite classifying space~$BG$
may be deduced from a statement that concerns only the large scale
geometry of~$G$.  In its most common formulation:

\begin{thm}[Descent]  \label{the:old_descent}
  Let~$G$ be a discrete group with finite classifying space~$BG$.  If
  the coarse Baum-Connes assembly map $\mu_X\colon \KX_*(G)\to\K_*(C^*
  (\abs{G}))$ is an isomorphism, then the assembly map $\mu\colon
  \K_*(BG)\to\K_*(\Cred(G))$ is injective.
\end{thm}

In this paper, we arrive at a new and somewhat strengthened version of
the Principle of Descent which deals with the existence of a
$\gamma$\nbd{}element.

We first recall how the notion of $\gamma$\nbd{}element arose in the
work of Kasparov in~\cite{Kasparov}.  Suppose that~$G$ is the
fundamental group of a compact, aspherical, even dimensional
spin$^{\mathrm{c}}$\brd{}manifold~$M$ and consider the class $\Dirac =
[\,\dirac\,]$ in $\K_0(M) \cong \KK^G (C_{0}(\tilde{M}),\C)$ of the
Dirac operator on~$M$.  If~$A$ is any $G$\nbd{}$C^*$\brd{}algebra,
taking the Kasparov product of~$\Dirac$ with~$1_A$, one obtains
a class $\Dirac\otimes 1_A\in \KK(C_0(\tilde{M})\otimes A\rcross G, A\rcross
G)$.  This class induces a map $\K_*\bigl(C_0(\tilde{M})\otimes A\rcross
G\bigr) \to \K_*(A \rcross G)$, which can be identified with the
Baum-Connes assembly map for~$G$ with coefficients in~$A$.  Under
certain geometric hypotheses, there exists
$\eta\in\KK^G(\C,C_0(\tilde{M}))$ such that $\Dirac
\otimes_\C \eta = 1_{C_0(\tilde{M})}$.  It follows that the
Baum-Connes assembly map is \emph{split injective}, and hence that the
Novikov conjecture holds for any group $G$ satisfying these geometric
hypotheses.  This method of verifying the Novikov conjecture has since
become known as the \emph{Dirac dual Dirac method}.

It is shown in~\cite{MeyerNest} that a substitute
$\Dirac\in\KK^G(\Proj,\C)$ for the class of the Dirac operator
in the above argument exists for any
locally compact group~$G$.  This class, or as we will refer to it, \emph{morphism}, 
is constructed using general
results on triangulated categories, and called a \emph{Dirac morphism}
for~$G$.  It is interpreted as a \emph{projective resolution} of~$\C$
in the category $\KK^G$, with respect to a certain localizing
subcategory.  For any $G$\nbd{}$C^*$\brd{}algebra~$A$, the morphism
$\Dirac\otimes 1_A$ induces a map $\K_*(\Proj \otimes A \rcross G) \to
\K_*(A\rcross G)$, and this map is naturally isomorphic to the
Baum-Connes assembly map.  This means that the left hand side
$\Ktop(G;A)$ of the Baum-Connes assembly map may be regarded as the
\emph{left derived functor} of the right hand side, $\K_*(A \rcross
G)$.  Projective resolutions are unique up to the obvious notion of
$\KK^G$\brd{}equivalance, so that~$\Dirac$ is canonically associated
to the group~$G$.

As a consequence, the Baum-Connes assembly map is an isomorphism with
arbitrary coefficients if~$\Dirac$ is invertible in $\KK^G(\Proj,\C)$,
and split injective with arbitrary coefficients if there exists
$\eta\in\KK^G(\C,\Proj)$ with $\Dirac \otimes_\C\eta=1_\Proj$.  Such
an~$\eta$ is called a \emph{dual Dirac morphism} and
$\gamma=\eta\otimes_\Proj \Dirac$ a \emph{$\gamma$\nbd{}element} for
$G$.  The dual Dirac morphism and the $\gamma$\nbd{}element are unique
if they exist.  If the Dirac dual Dirac method applies to~$G$ in the
sense of, say, \cite{Hig1}, then there exists a dual Dirac morphism in
the above sense.  The converse also holds if the domain~$\Proj$ of the
Dirac morphism can be taken to be proper.  This is satisfied for all
groups which we consider in this article. 

It is rather easy to see that a dual Dirac morphism exists if and only if
the Dirac morphism induces an isomorphism
$$
\Dirac\otimes{\cdot}\colon \KK^G (\C , \Proj) \to \KK^G (\Proj , \Proj)
$$
If~$G$ is discrete with a finite model for its classifying space~$BG$
then this is the case if and only if the Dirac morphism induces an
isomorphism
$$
\Dirac\otimes{\cdot}\colon \KK^G (\C,C_0(G)) \to \KK^G (\Proj,C_0(G)).
$$
We call this map the \emph{analytic co\brd{}assembly map}.  The main
observation of this article is that this map is a metric invariant
of~$G$.

To show this, we begin by identifying the domain $\KK^G(\C,C_0(G))$ of the above
map with
the reduced $\K$\nbd{}theory of a certain $C^*$\nbd{}algebra $\bv(G)$
which is a variant of the Higson corona of~$G$ and that only depends
on the structure of~$G$ as a coarse space and not on the algebraic
structure of~$G$ as a group.  We call $\bv(G)$ the \emph{stable Higson
corona} of~$G$.

The right hand side of the map above can be computed via a spectral
sequence with $E^2_{pq}=H^p(G,\Z G)$ for even~$q$ and $E^2_{pq}=0$ for
odd~$q$ and with $d^2_{pq}=0$.  It is known that the group $H^*(G,\Z
G)$ is isomorphic to the coarse cohomology $\HX^*(G)$ of~$G$ as
defined in~\cite{Roe} and in particular depends only on the structure
of~$G$ as a coarse space.  Suppose now that there exists a compact
model for the classifying space~$BG$.  This implies that~$G$ is
torsion free and that~$G$ is coarsely equivalent to~$\EG$.  Hence
$\HX^* (G) \cong \HX^*(\EG)$.  Since~$\EG$ is uniformly contractible,
$\HX^* (\EG) \cong H^*_{\mathrm{c}} (\EG)$, which is rationally
isomorphic to $\K^*(\EG)$.  These calculations suggest that
$\KK^G(\Proj,C_0(G))\cong \K^*(\EG)$, and this is indeed the case.

By
construction of the stable Higson corona there is a connecting map
$$
\mu^*\colon \tilde{\K}_{*+1}(\bv(G)) \overset{\cong}\to
\tilde{\K}_{*+1}(\bv(\EG)) \to \K^*(\EG).
$$
We call it the \emph{coarse co\brd{}assembly map for~$G$} due to its
duality, discussed in~\cite{EmersonMeyer}, with the ordinary coarse Baum-Connes assembly map appearing in
Theorem~\ref{the:old_descent}. To summarize, there are natural isomorphisms
for which the following diagram commutes:
$$
\xymatrix{
  {\KK^G_*(\C,C_0(G))} \ar[r]^{\Dirac\otimes{\cdot}} \ar[d]^{\cong} &
  {\KK^G _*(\Proj,C_0(G))} \ar[d]^{\cong} \\
  {\tilde{\K}_{*+1}(\bv(G))} \ar[r]^{\mu^*} &
  {\K^*(\EG).}
}
$$
The bottom map in this diagram is coarse, and the top determines
whether or not~$G$ has a dual Dirac morphism.  We have, therefore: 

\begin{thm}  \label{the:dual_Dirac_finite_BG}
  Let~$G$ be a discrete group with finite~$BG$.  Then~$G$ has a dual
  Dirac morphism if and only if the coarse co\brd{}assembly map $\mu^*
  \colon \tilde{\K}_{* +1}(\bv(G)) \to \K^*(\EG)$ is an isomorphism.
\end{thm}

When~$G$ does not possess a finite model for~$BG$ it becomes necessary
to coarsen the $\K$\nbd{}theory group appearing as the target of the
coarse co\brd{}assembly map $\mu^*$.  We define the
\emph{coarse $\K$\nbd{}theory of~$G$} by $\KX^* (G) \defeq \K_*\bigl(
\varprojlim C_0(P_d G ) \bigr)$, where $P_d G $ is the Rips complex
for $G$ of parameter $d$ and $\varprojlim C_0(P_d G )$ is the
$\sigma$\nbd{}$C^*$\brd{}algebra associated to the projective system
$C_0(P_{d+1}G) \to C_0(P_{d}G)$.  $\K$\nbd{}theory for
$\sigma$\nbd{}$C^*$\brd{}algebras is defined by N.~C.\ Phillips
in~\cite{Phillips}.  We are also obliged to introduce coefficients
into these groups and maps, but this is easily done, and in this way
one obtains a \emph{coarse co\brd{}assembly map with coefficients in a
$C^*$\nbd{}algebra~$D$},
$$
\mu^*_D \colon \tilde{\K}_{*+1}(\bv(G,D))\to\KX^*(G,D).
$$

\begin{thm}  \label{the:dual_Dirac_finite_dim_BG}
  Let~$G$ be a discrete, torsion free group with finite dimensional
  model for~$BG$.  Then~$G$ has a dual Dirac morphism if and only if
  the coarse co\brd{}assembly map~$\mu^*_D$ with coefficients in~$D$
  is an isomorphism for every $C^*$\nbd{}algebra~$D$.
\end{thm}

In fact, it suffices to assume that~$\mu^*_D$ is an isomorphism for
$D=C_0(\N)$.

\begin{cor}  \label{cor:gamma_coarse}
  Let~$G$ be a discrete, torsion free group with a finite dimensional
  model for~$BG$.  Then the existence or non\brd{}existence of a
  $\gamma$\nbd{}element for~$G$ is geometric, that is, it only
  depends on the large scale geometry of~$G$.  In particular, if two
  such groups $G$ and~$G'$ are coarsely equivalent, then~$G$ has a
  $\gamma$\nbd{}element if and only if~$G'$ does.
\end{cor}

The hypothesis of finite dimensional~$BG$ is needed because we use a
concrete model for the Dirac morphism constructed by G.\ Kasparov and G.\
Skandalis in~\cite{KasparovSkandalis}.  It would be more in the spirit
of~\cite{MeyerNest} to work with the abstract definition of the Dirac
morphism.  However, issues with countably infinite direct sums make such
a line of argument technically difficult.  The problem only occurs in
the passage from the analytic co\brd{}assembly map to the existence of
a dual Dirac morphism, which involves a Mayer-Vietoris argument (see
Section~\ref{sec:geo_finite}).

We also investigate the case in which~$G$ has torsion.  To do so, we
work equivariantly with respect to finite subgroups of~$G$.  Any such
subgroup $H\subseteq G$ acts on $\bv(G)$, so that we can form the
crossed product $\bv(G)\cross H$.  More generally, if~$D$ is an
$H$\nbd{}$C^*$\brd{}algebra, then~$H$ acts on $\bv(G,D)$.  We
construct an \emph{$H$\nbd{}equivariant coarse co\brd{}assembly map
with coefficients in $D$}
$$
\mu^*_{D,H}\colon \tilde{\K}_{*+1}(\bv(G,D)\cross H) \to \KX^*_H(G,D).
$$
This map depends only on the $H$\nbd{}equivariant coarse equivalence
class of~$G$.  If~$G$ possesses a $\gamma$\nbd{}element, then
$\mu^*_{D,H}$ is an isomorphism for all $D$ and~$H$.  With our current
methods, we can prove the converse under some hypotheses on~$G$:

\begin{thm}  \label{the:dual_Dirac_finite_dim_EG}
  Let~$G$ be a discrete group with a finite dimensional model
  for~$\EG$, and assume that $G$ has only finitely many conjugacy
  classes of finite subgroups.  (This occurs for instance if~$G$ has a
  $G$\nbd{}finite model for~$\EG$.)  Then~$G$ has a dual Dirac
  morphism if and only if the $H$\nbd{}equivariant coarse
  co\brd{}assembly map with coefficients in~$D$ is an isomorphism for
  every finite subgroup~$H$ and every $H$\nbd{}$C^*$\brd{}algebra~$D$.
\end{thm}

Once again, it suffices to require an isomorphism for $D=C_0(\N)$.

Let~$G$ be the fundamental group of a compact, aspherical
manifold~$M$.  Then the $\K$\nbd{}theory of~$\tilde{M}$ is generated
by the Bott class $\beta\in\K^{-n}(\tilde{M})$.
Theorem~\ref{the:dual_Dirac_finite_dim_EG} implies that~$G$ has a dual
Dirac morphism if and only if there is a unique~$\theta$ in the
$\K$\nbd{}theory of $\bv(G)$ with $\partial(\theta)=\beta$.  If this
is the case, then the dual Dirac morphism can be manufactured in
finitely many steps from~$\theta$, each of them using a Mayer-Vietoris
argument. The condition that $\partial (\theta ) = \beta$ Poincar\'e dualizes 
to ensure that $\Dirac \otimes \eta = 1$. Similar remarks hold for 
general (torsion free) discrete groups $G$ with finite-dimensional $\EG$. Thus, a dual-Dirac
morphism for a group $G$ must, if it exists, arise from $\K$-theory classes
for the stable Higson corona of
$G$.

The coarse co\brd{}assembly map for a general coarse space~$X$ is
studied in more detail in~\cite{EmersonMeyer}, where we give several
cases in which it is an isomorphism.  This is the case for scaleable
spaces in the sense of~\cite{HigsonRoe} and for spaces which uniformly
embed in a Hilbert space.  The first result implies that a torsion
free group that is coarsely equivalent to a uniformly contractible,
scaleable space, has a $\gamma$\nbd{}element.  Such a
$\gamma$\nbd{}element \emph{a priori} does not arise from any of the
usual methods of constructing them (for example from Lipschitz maps to
Euclidean space).

\section{Projective resolutions, Dirac and dual Dirac morphisms}
\label{sec:resolutions_Dirac_dual_Dirac}

We shall use some ideas of~\cite{MeyerNest}, in which the role of the
Dirac operator is centralized in the construction of the Baum-Connes
assembly map as a derived functor.  We begin by recalling these
results.

Let~$G$ be a locally compact group, $H$ a compact subgroup of~$G$,
and~$D$ an $H$\nbd{}$C^*$\brd{}algebra.  One has two functors: the
\emph{restriction functor} $\Res_G^H\colon \KK^G\to\KK^H$, whose
definition is obvious, and the \emph{induction functor}
$\Ind_H^G\colon\KK^H\to\KK^G$.  The latter is defined on objects, that
is, $H$\nbd{}$C^*$\brd{}algebras, by setting
$$
\Ind_H^G D \defeq
\{f\in C_0(G,D) \mid
  \text{$\alpha_h\bigl(f(g)\bigr)=f(gh)$ for all $h\in H$, $g\in G$}
\},
$$
with $G$\nbd{}action $(gf)(g') = f(g^{-1}g')$.  Similarly one defines
$\Ind_H^G(\Hilm)$ for $H$\nbd{}equivariant Hilbert (bi)modules.  We
call an object of $\KK^G$ \emph{compactly induced} if it is
$\KK^G$\brd{}equivalent to a $G$\nbd{}$C^*$\brd{}algebra of the form
$\Ind_H^G D$.  Let $\CI\subseteq\KK^G$ denote the class of all
compactly induced objects, and let $\gen{\CI}\subseteq\KK^G$ be the
localizing subcategory generated by~$\CI$.  This is, by definition,
the smallest full subcategory of~$\KK^G$ satisfying:
\begin{enumerate}[(1)]
\item $\gen{\CI}$ contains~$\CI$;

\item $\gen{\CI}$ is triangulated, that is, closed under suspensions
  and under extensions with an equivariant, completely positive
  section;

\item $\gen{\CI}$ is closed under countable direct sums.

\end{enumerate}
The subcategory~$\gen{\CI}$ contains all proper
$G$\nbd{}$C^*$\brd{}algebras (see~\cite{MeyerNest}).  The Baum-Connes
assembly map is an isomorphism for coefficients in~$\gen{\CI}$.

An element $f\in\KK^G(A,B)$ is called a \emph{weak equivalence} if
$\Res_G^H (f)$ is invertible in $\KK^H(A,B)$ for all compact subgroups
$H\subseteq G$.  An object $A\in\KK^G$ is called \emph{weakly
contractible} if $\Res_G^H(A)\cong0$.  Equivalently, the zero map
$0\to A$ is a weak equivalence.  We let $\CC\subseteq\KK^G$ be the
full subcategory of weakly contractible objects.

The idea of \cite{MeyerNest} is to mimic the construction of derived
functors and categories in homological algebra using the above
definitions.  This leads to a sort of dictionary in which the weakly
contractible objects play the role of the exact chain complexes, the
weak equivalences play the role of the quasi-isomorphisms, and the
objects of~$\gen{\CI}$ play the role of the projective chain
complexes.  Thus the analogue of a projective resolution of $A\in\KK^G$
is a weak equivalence $P\to A$ with $P\in\gen{\CI}$.  The case $A=\C$
is especially important:

\begin{defn}  \label{def:Dirac}
  A \emph{Dirac morphism} for~$G$ is a weak equivalence
  $\Dirac\in\KK^G(\Proj,\C)$ with $\Proj\in\gen{\CI}$.
\end{defn}

In \cite{MeyerNest}, the following facts are established:

\begin{thm}  \label{the:MeyerNest}
  Let~$G$ be a locally compact group.
  \begin{enumerate}[(1)]
  \item A Dirac morphism for~$G$ exists and is unique up to
    $\KK^G$\brd{}equivalence.
    
  \item The Baum-Connes assembly map with coefficients in~$A$ is
    naturally isomorphic to the map
    $$
    \Dirac_*\colon \K_*((\Proj\otimes A)\rcross G)\to\K_*(A\rcross G),
    $$
    induced by a Dirac morphism $\Dirac\in\KK^G(\Proj,\C)$.

  \item Let~$\EG$ be a locally compact model for the universal proper
    $G$\nbd{}space and let $\Dirac\in\KK^G(\Proj,\C)$ be a Dirac
    morphism.  For every pair of $G$\nbd{}$C^*$\brd{}algebras $A$
    and~$B$ there is a natural isomorphism
    $$
    \delta_{AB}\colon \KK^G(\Proj\otimes A,B) \cong \RKK^G(\EG;A,B)
    $$
    such that the following diagram commutes:
    \begin{displaymath}
      \xymatrix@C=-1em{
        {\KK^G(A,B)} \ar[rr]^{\delta_{AB}}_{\cong} & &
        {\RKK^G(\EG;A,B)} \\
        & {\KK^G(A,B),}  \ar[ur]_{p_\EG^*}
        \ar[ul]^{\Dirac\otimes{\cdot}} &
      }
    \end{displaymath}
    where $p_\EG^*\colon \KK^G(A,B)\to\RKK^G (\EG;A,B)$ is induced
    by the map $\EG\to\star$, where~$\star$ is a point.
  \end{enumerate}
\end{thm}

Statement~(1) asserts that projective resolutions exist for all
$A\in\KK^G$.  Statement~(2) implies that the functor $A \mapsto \Ktop
(G,A)$ is the left derived functor of the functor $A \to \K(A \rcross
G)$ and statement~(3) means that the category $\RKK^G(\EG)$ is the
derived category of $\KK^G$ with respect to the weak equivalences.
The third statement also allows us to give the following
characterization of the Dirac morphism, which will be of use to us.

\begin{lemma}  \label{lem:poincareduality}
  Let~$G$ be a locally compact group, let~$A$ be a
  $G$\nbd{}$C^*$\brd{}algebra and $d\in\KK^G(A,\C)$.  Then~$d$ is a
  Dirac morphism for~$G$ if and only if for every $G$\nbd{}$C^*$\brd{}algebra $B$ 
there is a natural isomorphism
  $$
  \delta_B\colon \KK^G(A,B) \cong \RKK^G(\EG;\C,B)
  $$such that the following diagram commutes:
  \begin{equation}\label{poincaredualitydiagram}
    \xymatrix@C=-1em{
      {\KK^G(A,B)} \ar[rr]^{\delta_B}_{\cong} & &
      {\RKK^G(\EG;\C,B)} \\
      & {\KK^G(\C,B).}  \ar[ur]_{p_\EG^*} \ar[ul]^{d\otimes{\cdot}} &
    }
  \end{equation}
\end{lemma}

\begin{proof}
  It is shown in~\cite{MeyerNest} that a Dirac morphism
  $\Dirac\in\KK^G(\Proj,\C)$ has these properties (compare~(3) in
  Theorem~\ref{the:MeyerNest}).  Conversely, the hypotheses on~$d$
  determine the functor $B\mapsto \KK^G(A,B)$ and the natural
  transformation $d\otimes{\cdot}\colon \KK^G(\C,B)\to\KK^G(A,B)$
  uniquely.  By the Yoneda Lemma, this implies that $A\cong\Proj$ and
  that~$d$ corresponds to~$\Dirac$ under this isomorphism.
\end{proof}

The isomorphism~\eqref{poincaredualitydiagram} holds for the case of
the classical Dirac operator on a manifold (see
\cite{Kasparov}*{Theorem 4.9}) and for the Kasparov-Skandalis Dirac
morphism associated to a finite dimensional simplicial complex (see
\cite{KasparovSkandalis}*{Theorem 6.5}).  Thus we have the following.

\begin{cor}  \label{realizationsofthedirac}
  Let~$G$ be a locally compact group.
  \begin{enumerate}[A.]
  \item If~$\EG$ can be realized by a complete Riemmanian manifold on
    which~$G$ acts properly and isometrically, then the class $[\Dirac_X]
    \in\KK^G(C_\tau(X),\C)$ of~\cite{Kasparov}*{Definition 4.2} is a
    Dirac morphism for~$G$.
    
  \item If~$\EG$ can be realized by a finite dimensional simplicial
    complex~$X$ on which~$G$ acts simplicially, then the Dirac
    morphism for~$G$ may be identified with the class $[\Dirac_X]\in
    \KK^G(\Proj_X,\C)$ of \cite{KasparovSkandalis}*{Definition 1.3}.

  \end{enumerate}
\end{cor}

\begin{note}
  In Corollary~\ref{realizationsofthedirac}.B we have changed notation
  from \cite{KasparovSkandalis}, denoting by~$\Proj_X$ what they have
  denoted by~$\mathcal{A}_X$.
\end{note}

\begin{rmk}
  Formally, the source~$\Proj$ of the Dirac morphism must be an
  ungraded $G$\nbd{}$C^*$\brd{}algebra because the
  $\KK^G$\brd{}category of graded $C^*$\brd{}algebras is not
  triangulated.  However, it is certainly permissible to use a graded
  $G$\nbd{}$C^*$\brd{}algebra that is $\KK^G$\brd{}equivalent to an
  ungraded one.  It is well-known that the $G$\nbd{}$C^*$\brd{}algebra
  $C_\tau(X)$ in Corollary~\ref{realizationsofthedirac}.A is
  $\KK^G$\brd{}equivalent to $C_0(T^* M)$.  A similar ungraded model
  for~$\Proj_X$ is constructed in~\cite{KasparovSkandalis}.
  Therefore, we may ignore this technical issue in the following.
\end{rmk}

\begin{defn}  \label{def:dual_Dirac}
  Let~$G$ be a locally compact group and let
  $\Dirac\in\KK^G(\Proj,\C)$.  A \emph{dual Dirac morphism} for~$G$ is
  an element $\eta\in\KK^G(\C,\Proj)$ such that $\Dirac\otimes_\C\eta
  = 1_\Proj$.  The composition $\gamma = \eta \otimes_\Proj \Dirac$ is
  called a \emph{$\gamma$\nbd{}element} for~$G$.
\end{defn} 

By definition, $G$ has a $\gamma$\nbd{}element if and only if~$G$ has
a dual Dirac morphism.  Moreover, if a dual Dirac morphism exists,
then it is unique.  Consequently the same is true of
$\gamma$\nbd{}elements.

\begin{rmk}  \label{splitting_KKG}
  It is shown in~\cite{MeyerNest} that a dual Dirac morphism exists if
  and only if $\KK^G$ is isomorphic as a triangulated category to the
  product of the subcategory~$\gen{\CI}$ and the subcategory of weakly
  contractible objects~$\CC$.  For $A\in\KK^G$, let $\gamma_A\defeq
  \gamma\otimes_\C 1_A\in\KK^G(A,A)$.  Then $A\in\gen{\CI}$ if and
  only if $\gamma_A=1_A$, and $A\in\CC$ if and only if $\gamma_A=0$.
  Thus the existence of a dual Dirac morphism implies that there is no
  interaction between $\gen{\CI}$ and~$\CC$, a situation which we note has no
  analogue in homological algebra.
\end{rmk}

The following is also shown in~\cite{MeyerNest} and follows easily
from Remark~\ref{splitting_KKG}.

\begin{prop}  \label{strengthening}
  If a dual Dirac morphism exists, then the map
  $$
  p_\EG^*\colon \KK^G(A,B)\to\RKK^G(\EG;A,B)
  $$
  is an isomorphism for all $A\in\KK^G$, $B\in\gen{\CI}$.
\end{prop}

For convenience, we introduce the following definition.

\begin{defn}
  The \emph{analytic co\brd{}assembly map with coefficients in $A$} is
  the map
  \begin{equation}  \label{analyticcoassemblymap}
    p_\EG^*\colon \KK^G(\C,A)\to\RKK^G(\EG;\C,A),
  \end{equation}
  where $p_\EG\colon \EG\to\star$ is the constant map from~$\EG$ to a point.
\end{defn}

\begin{prop}  \label{pro:dual_Dirac_isomorphism}
  Let~$G$ be a locally compact group and $\Dirac\in\KK^G(\Proj,\C)$
  a Dirac morphism for~$G$.  Then~$G$ possesses a dual Dirac
  morphism if and only if the analytic co\brd{}assembly map with
  coefficients in~$\Proj$
  $$
  p_\EG^*\colon \KK^G(\C,\Proj)\to\RKK^G(\EG;\C,\Proj)
  $$
  is an isomorphism.  
\end{prop}

\begin{proof}
  By statement~(3) of Theorem~\ref{the:MeyerNest}, this map is an
  isomorphism if and only if the map
  $$
  \Dirac^*\colon \KK^G(\C,\Proj)\to\KK^G(\Proj,\Proj)
  $$
  is an isomorphism.  If this is the case, then the inverse
  image~$\eta$ of $1_\Proj\in\KK^G(\Proj,\Proj)$ is a dual Dirac
  morphism.  The converse implication is contained in
  Proposition~\ref{strengthening}, since by definition $\Proj \in \gen{\CI}$.

\end{proof}

\begin{rmk}
  Another definition of a $\gamma$\nbd{}element that is used
  frequently is the following (see~\cite{Hig1}).  An element
  $\gamma\in\KK^G(\C,\C)$ is called a $\gamma$\nbd{}element if there
  is a factorization $\gamma = \alpha\otimes_P\beta$, where~$P$ is a
  proper $G$\nbd{}$C^*$\brd{}algebra, $\alpha\in\KK^G(\C,P)$,
  $\beta\in\KK^G(P,\C)$, and $\beta\otimes_\C\alpha=1_P$.  If such
  $P$, $\alpha$ and $\beta$ exist, then~$\alpha$ is necessarily a weak
  equivalence and hence a model for the Dirac morphism of~$G$.
  Hence~$\beta$ is a dual Dirac morphism in our sense, and likewise
  for~$\gamma$.  We do not, however, know whether the converse is
  true.  That is, we do not know whether the source $\Proj$ of the
  Dirac morphism can be realized in general by a proper
  $G$\nbd{}$C^*$\brd{}algebra.  Hence our definition is \emph{a
  priori} less restrictive than that of~\cite{Hig1}.
\end{rmk}

\section{Geometrically finite groups}
\label{sec:geo_finite}

Propositions \ref{strengthening} and~\ref{pro:dual_Dirac_isomorphism}
show that a dual Dirac morphism exists if and only if the analytic
co\brd{}assembly map~\eqref{analyticcoassemblymap} is an isomorphism
for all coefficients in~$\gen{\CI}$ or, equivalently, for the fixed
coefficient~$\Proj$, where $\Dirac\in\KK^G(\Proj,\C)$ is a Dirac
morphism.  Since the category $\gen{\CI}$ is by definition generated
by~$\CI$, it seems plausible that it suffices to check the isomorphism
on objects of~$\CI$.  However, since we know nothing about the
behavior of $\KK^G$ under infinite direct sums in the second variable,
there is a difficulty in passing from~$\CI$ to $\gen{\CI}$.  To avoid
this difficulty, we restrict attention to groups that satisfy some
finiteness hypotheses that allow us to construct the domain~$\Proj$ of
the Dirac morphism from compactly induced $G$\nbd{}$C^*$\brd{}algebras
without using infinite direct sums.  This technique already covers
several cases of interest.

\begin{defn} \label{definitionofgeometricallyfinite}
  Let~$G$ be a discrete, countable group.
  \begin{enumerate}[A.]
  \item We say that~$G$ is \emph{geometrically finite} if~$\EG$ can be
    realized as a finite dimensional simplicial complex. 

  \item We say that~$G$ is \emph{strongly geometrically finite} if it is
    geometrically finite and, in addition, has at most finitely many
    conjugacy classes of finite subgroups.
  \end{enumerate}
\end{defn}

\begin{rmk}
  Clearly, if~$G$ has a $G$\nbd{}finite model for~$\EG$, then~$G$ is
  strongly geometrically finite, but the condition of strong geometric
  finiteness is obviously much weaker.
\end{rmk}

\begin{thm}  \label{bigtheorem2}
  Let~$G$ be a discrete group.
  \begin{enumerate}[A.]
  \item If~$G$ is strongly geometrically finite, then~$G$ has a
    $\gamma$\nbd{}element if and only if the analytic co\brd{}assembly
    map~\eqref{analyticcoassemblymap} with coefficients in $\Ind_H^G D$ is an
    isomorphism for every finite subgroup $H\subseteq G$ and every
    $H$\nbd{}$C^*$\brd{}algebra~$D$.  

  \item If~$G$ has a $G$\nbd{}finite model for~$\EG$, then~$G$ has a
    $\gamma$\nbd{}element if and only if the analytic co\brd{}assembly
    map~\eqref{analyticcoassemblymap} with coefficients in $C_0(G/H)$ is an
    isomorphism for every finite subgroup $H\subseteq G$.

  \end{enumerate}
\end{thm}

\begin{rmk}
  The proof shows that in case~A a $\gamma$\nbd{}element already
  exists if the analytic co\brd{}assembly map is an isomorphism for
  $C_0(G/H\times\N)\cong \Ind_H^G C_0(\N)$ for any~$H$.  Thus we only
  need one very simple coefficient algebra.  Since~$\C$ is a direct
  summand of $C_0(\N)$, the analytic co\brd{}assembly map for
  $C_0(G/H)$ is a direct summand of the analytic co\brd{}assembly map
  for $C_0(G/H\times \N)$.  Therefore, if we have an isomorphism for
  $C_0(G/H\times\N)$, then we also have an isomorphism for~$\C$.
\end{rmk}

\begin{proof}
  We have already observed in
  Proposition~\ref{pro:dual_Dirac_isomorphism} that the existence of a
  $\gamma$\nbd{}element implies that the analytic co\brd{}assembly map
  is an isomorphism for all coefficients in~$\gen{\CI}$, without any
  hypothesis on the group~$G$.  We have to prove the converse.
  Let~$G$ be geometrically finite, let~$X$ be a finite dimensional
  simplicial complex realizing~$\EG$, and let
  $[\Dirac_X]\in\KK^G(\Proj_X,\C)$ be the Kasparov-Skandalis
  realization of the Dirac morphism for~$G$.  By
  Proposition~\ref{pro:dual_Dirac_isomorphism} and
  Corollary~\ref{realizationsofthedirac}, $G$ possesses a
  $\gamma$\nbd{}element if and only if the analytic co\brd{}assembly
  map with coefficients in~$\Proj_X$ is an isomorphism.  Recall
  from~\cite{KasparovSkandalis} that the skeletal filtration of~$X$
  gives rise to a filtration of~$\Proj_X$ by ideals
  $$
  0 = \Proj^{(-1)}_X \subset
  \Proj^{(0)}_X \subset
  \Proj^{(1)}_X \subset
  \Proj^{(2)}_X \subset \cdots
  \Proj_X^{(n)} =
  \Proj_X.
  $$
  Since the resulting extensions
  $$
  0 \to
  \Proj_X^{(k-1)} \to
  \Proj_X^{(k)} \to
  \Proj_X^{(k)}/\Proj_X^{(k-1)} \to
  0
  $$
  have $G$\nbd{}equivariant completely positive sections, $\Proj_X$
  lies in the triangulated subcategory of $\KK^G$ that is generated by
  the subquotients $\Proj_X^{(k)}/\Proj_X^{(k-1)}$.  The latter are
  $\KK^G$\nbd{}equivalent to $C_0(X^{(k)})$, where~$X^{(k)}$ denotes
  the set of $k$\nbd{}cells of~$X$, viewed as a discrete
  $G$\nbd{}space.  It follows that the analytic co\brd{}assembly map
  is an isomorphism with coefficients in~$\Proj_X$ if it is an
  isomorphism with coefficients in $C_0(X^{(k)})$ for all $k\in\N$.
  Another way of expressing this is as follows.  The extensions above
  give rise to commutative diagrams with exact rows:
  $$
  \xymatrix@C=0.6em{
    {} \ar[r] &
    {\KK^G_*(\C,\Proj_X^{(k-1)})} \ar[d] \ar[r] &
    {\KK^G_*(\C,\Proj_X^{(k)})} \ar[r] \ar[d] &
    {\KK^G_*(\C,\Proj_X^{(k)} / \Proj_X^{(k-1)})} \ar[d] \ar[r] &
    {} \\
    {} \ar[r] &
    {\RKK^G_*(\EG;\C,\Proj_X^{(k-1)})}  \ar[r] &
    {\RKK^G_*(\EG;\C,\Proj_X^{(k)})} \ar[r] &
    {\RKK^G_*(\EG;\C,\Proj_X^{(k)}/\Proj_X^{(k-1)})} \ar[r] &
    {,}
  }
  $$
  where the vertical maps are induced by product with the Dirac
  morphism $[\Dirac_{X}]$.  If we apply the Five Lemma to these
  diagrams, we obtain by induction on~$k$ that the analytic
  co\brd{}assembly map is an isomorphism with coefficients
  $\Proj_X^{(k)}$ for $k=-1,\dots,n$ provided it is an isomorphism for
  the subquotients.  By the identification of the subquotients above,
  this will be the case if the analytic co\brd{}assembly map with
  coefficients $C_0(X^{(k)})$ is an isomorphism for all $k$.
  
  Each of the discrete, proper $G$\nbd{}spaces~$X^{(k)}$ is
  $G$\nbd{}isomorphic to a countably infinite disjoint union of
  homogeneous spaces $G/H$ for finite subgroups $H\subseteq G$.
  Since~$G$ has only finitely many conjugacy classes of finite
  subgroups, at most finitely many of these homogeneous spaces are
  non-isomorphic as $G$\nbd{}spaces.  Thus $C_0(X^{(k)})$ is
  $G$\nbd{}isomorphic to a \emph{finite} direct sum of
  $G$\nbd{}$C^*$\brd{}algebras of the form $C_0(I\times G/H)$ with
  finite subgroups $H\subseteq G$ and a discrete set~$I$ with trivial
  action of~$G$.  In this way we are reduced to verifying isomorphism
  of the analytic co\brd{}assembly map with coefficients $C_0(G/H
  \times I)$ for an at most countable set~$I$.  In case~B, only finite sets $I$ 
are involved, so that in this case we only need the analytic co\brd{}assembly
  map to be an isomorphism with coefficients $C_0(G/H)$.  If the
  set~$I$ is infinite, we can identify it with the fixed, countable
  set~$\N$ with trivial $G$\nbd{}action, so that we only require the
  analytic co\brd{}assembly map to be an isomorphism with coefficients
  $C_0(G/H\times\N)$.
\end{proof}

\section{The stable Higson corona}
\label{sec:stable_corona}

For basic matters on coarse structures we refer the reader
to~\cite{HigsonRoe}.  Such a structure on a locally compact
topological space~$X$ is given by a collection of \emph{entourages}
$E\subset X\times X$ satisfying various axioms.  A subset $B\subset X$
is called \emph{bounded} (with respect to the coarse structure) if
$B\times B$ is an entourage.  The coarse structure is called
\emph{proper} if every bounded subset is compact, and \emph{separable}
if~$X$ has a countable, uniformly bounded open cover.  Finally, we say
that the coarse structure is \emph{unital} if the diagonal is an
entourage.

An example of a coarse structure satisfying all these conditions
arises when~$X$ is a separable metric space of uniformly bounded
geometry, and the entourages are given by sets of the form
$$
E=\{(x,y)\in X\times X \mid d(x,y) \le R\}.
$$
We call a coarse structure metrizable if it arises from a metric in
this fashion.  Any proper, unital, separable coarse structure on a
second countable space~$X$ is metrizable.  In this paper we are
exclusively interested in this case.

\emph{Convention: By ``coarse space'' we always mean a second countable,
  locally compact topological space equipped with a separable, proper,
  and unital coarse structure.}

A (not necessarily continuous) map between coarse spaces is a
\emph{coarse map} if it maps entourages in $X$ to entourages in $Y$
and is \emph{proper} in the sense that inverse images of bounded
subsets are bounded.  Two maps $\phi$ and~$\psi$ from~$X$ to~$Y$ are
called \emph{close} if $\{(\phi(x),\psi(x)) \mid x\in X \}$ is an
entourage in $Y\times Y$; that is, if $\phi\times\psi$ maps the
diagonal in $X\times X$ to an entourage in $Y\times Y$.  Finally, we
say that two coarse spaces $X$ and~$Y$ are \emph{coarsely equivalent}
if there exist coarse maps $X\to Y$ and $Y\to X$ whose compositions
are close to the identity maps on $X$ and~$Y$, respectively.  Finally,
if~$X$ is a coarse space and~$G$ is a group acting by homeomorphisms
of $X$, then a coarse structure is $G$\nbd{}invariant if any entourage
is contained in a $G$\nbd{}invariant entourage.  We also say that $G$
acts \emph{isometrically} on the coarse space~$X$.

Every discrete and countable group~$G$ may be regarded as a coarse
space.  Its coarse structure is defined by the entourages
$$
E = \{(g_1,g_2) \mid g_1g_2^{-1} \in F \},
$$
where~$F$ ranges over the finite subsets of~$G$.  This coarse
structure is manifestly \emph{right} translation invariant and is the
unique coarse structure on~$G$ with this property.  We always
equip~$G$ with this coarse structure.  The metrizability of this
coarse structure means that there is a proper length function~$l$
on~$G$ such that the associated right invariant metric $d(g,h) =
l(gh^{-1})$ induces the coarse structure.  It is clear that any two
such $l$ induce the same coarse structure; thus the length function
$l$ is not unique, but the coarse structure it generates is.

Let $f\colon X \to E$ be a function from a coarse space $X$ into a
Banach space $E$. Define a function $\nabla f\colon X\times X \to\R_+$
by $\nabla f(x,y)=\norm{f(x)-f(y)}$.

\begin{defn}  \label{vanishingvariation}
  Let~$X$ be a coarse space and let~$f$ be a continuous function
  on~$X$ with values in a Banach space $E$.  Then we say that~$f$ has
  \emph{vanishing variation} if the restriction of $\nabla f$ to the
  closure of any entourage vanishes at infinity.
\end{defn}

Let~$\Comp$ be the $C^*$\nbd{}algebra of compact operators on
$\ell^2(\N)$.

\begin{defn}
  Let~$X$ be a coarse space and let~$D$ be a $C^*$\nbd{}algebra.  Then
  $\vv(X,D)$ shall denote the $C^*$\nbd{}algebra of bounded,
  continuous functions $X\to D\otimes\Comp$ of vanishing variation.
\end{defn}

The central definition of this paper is the following.

\begin{defn}
  Let~$X$ be a coarse space.  Then the \emph{stable Higson corona $\bv
  (X,D)$ of $X$ with coefficients in $D$} is the $C^*$\nbd{}algebra
  $$
  \bv (X,D) = \vv(G,D) / C_0(X,D\otimes\Comp).
  $$
\end{defn}

\begin{rmk}  \label{rem:corona_coarse}
  It is not difficult to show that the isomorphism class of the algebra
  $\bv(X,D)$ only depends on the coarse equivalence class of~$X$.  In
  particular, if one replaces~$X$ by a discrete subset~$X'$ such that
  the inclusion $X'\to X$ is a coarse equivalence, then the
  restriction map $\bv(X,D)\to\bv(X',D)$ is a $*$\nbd{}isomorphism.
  More generally, a coarse map $X\to Y$ induces a canonical
  $*$\brd{}homomorphism $\bv(Y,D)\to\bv(X,D)$, and two close maps
  induce the same homomorphism.  Hence $X\mapsto\bv(X,D)$ is
  functorial from the category of coarse spaces and coarse maps to the
  category of $C^*$\nbd{}algebras and $*$\brd{}homomorphisms.
  See~\cite{EmersonMeyer} for details, or~\cite{Roe} for the analogous
  assertions for the Higson corona.
\end{rmk}

It is convenient for this article to replace $\bv(X,D)$ and $\vv(X,D)$
by the following variants.  We denote the multiplier algebra of a
$C^*$\nbd{}algebra~$D$ by $\Mult(D)$ and the stable multiplier algebra
by $\Mult^s (D) \defeq\Mult (D \otimes \Comp)$. 

\begin{defn}
  Let~$X$ be a coarse space.  Then $\vvr(X,D)$ shall denote the
  $C^*$\nbd{}algebra of bounded, continuous functions of vanishing
  variation $f\colon X\to \Mult^s(D)$ such that
  $f(x)-f(y)\in D\otimes\Comp$ for all $x,y\in X$.  The \emph{reduced
  stable Higson corona of $X$ with coefficients in $D$, $\bvr(X,D)$}, is
  the quotient
  $$
  \bvr(X,D) = \vvr(X,D) / C_0(X,D\otimes\Comp).
  $$
\end{defn} 

The reason for the terminology is the following.  Let~$i$ be the
inclusion $D\otimes\Comp\to\vv(X,D)$ that sends elements of
$D\otimes\Comp$ to constant functions on~$X$.  We obtain an induced
map $i_*\colon\K(D)\cong\K(D\otimes\Comp)\to \K(\vv(X,D))$.  Let
$\pi\colon\vv(X,D)\to\bv(X,D)$ be the quotient map.

\begin{defn}
  We define the reduced $\K$\nbd{}theories of $\bv(X,D)$ and $\vv(X,D)$ by
  \begin{align*}
    \tilde{\K}_*(\bv(X,D)) &= 
    \K_*(\bv(X,D)) / (\pi \circ i)_* \K_*(D),
    \\
    \tilde{\K}_*(\vv(X,D)) & = 
    \K_*(\vv(X,D))/ i_* \K_*(D).
  \end{align*}
\end{defn}

\begin{lemma}[see~\cite{EmersonMeyer}]
  There are natural isomorphisms
  $$
  \tilde{\K}_*(\bv(X,D)) \cong \K_*(\bvr(X,D)),
  \qquad
  \tilde{\K}_*(\vv(X,D)) \cong \K_*(\vvr(X,D)).
  $$
\end{lemma}

From now on, we work exclusively with the reduced algebras $\vvr(X,D)$
and $\bvr(X,D)$.  However, it is the non-reduced algebras $\bv(X,D)$
and $\vv(X,D)$ which are most obviously connected to classical
constructions in coarse geometry, in particular to the Higson corona
construction (see \cites{HigsonRoe, Roe} for its definition.)  As
in Remark~\ref{rem:corona_coarse}, one can easily check the following
lemma.  Details can be found in~\cite{EmersonMeyer}.

\begin{lemma}
  The assignment $X \mapsto \bvr(X,D)$ is functorial from the category
  of coarse spaces and coarse maps to the category of
  $C^*$\nbd{}algebras and $C^*$\nbd{}algebra homomorphisms.  That is,
  a coarse map $X\to Y$ induces a canonical $*$\nbd{}homomorphism
  $\bvr(Y,D)\to\bvr(X,D)$, and close maps induce the same
  $*$\nbd{}homomorphism.
\end{lemma}

Let~$H$ be a finite group acting continuously on the coarse space~$X$
and preserving the coarse structure.  For compatibility with our later
arguments, we let~$H$ act on the right.  Suppose that~$H$ also acts on
the coefficient algebra~$D$ on the left via a homomorphism
$H\to\Aut(D)$, $h\mapsto\alpha_h$.  Then we obtain an action of~$H$ on
the algebra $\vvr(X,D)$ by $(h\cdot f)(x) = \alpha_h f(xh)$.  This
restricts to an action of~$H$ on the ideal $C_0(X,D\otimes\Comp)$ and
so descends to an action of~$H$ on $\bvr(X,D)$.  If~$H$ also acts on a
coarse space~$Y$, and $\phi\colon X\to Y$ is a coarse, proper,
$H$\nbd{}equivariant map, then the induced $*$\nbd{}homomorphism
$\bvr(Y,D)\to\bvr(X,D)$ is $H$\nbd{}equivariant, so that we obtain a
$*$\nbd{}homomorphism
$$
\bvr(X,D)\cross H\to\bvr(Y,D)\cross H.
$$
Since~$H$ is finite the type of cross product we use is immaterial.

By construction, there is an exact sequence
$$
0 \to
C_0(X,D\otimes\Comp) \to
\vvr(X,D) \to
\bvr(X,D) \to
0.
$$
If~$H$ is a finite group acting on~$X$ and preserving the coarse
structure, then we obtain an exact sequence
$$
0 \to
C_0(X,D\otimes\Comp)\cross H \to
\vvr(X,D)\cross H \to
\bvr(X,D)\cross H \to
0.
$$
This exact sequence is the origin of the \emph{coarse co\brd{}assembly
map} in Section~\ref{sec:coarse_coassembly}.

\section{Calculation of $\KK^G(\C,\Ind_H^G D)$}
\label{sec:aco_source}

Fix a discrete group~$G$, let~$H$ be a finite subgroup, and let~$D$ be
an $H$\nbd{}$C^*$\brd{}algebra.  We are going to identify the group
$\KK^G(\C,\Ind_H^G D)$ that appears in connection with
$\gamma$\nbd{}elements (Theorem~\ref{bigtheorem2}) with the
$\K$\nbd{}theory of the $C^*$\nbd{}algebra $\bvr(G,D)\cross H$
described in the previous section.  We let~$H$ act on~$G$ by right
translations.  This action preserves the coarse structure and hence
induces an action on the $C^*$\nbd{}algebra $\bvr(G,D)$.  We also let
$D_H\defeq D\otimes\Comp(\ell^2 H)$, equipped with the usual action
of~$H$, and denote the $H$\brd{}fixed point subalgebra of
$\bvr(G,D_H)$ by $\bvr(G,D_H)^H$.

\begin{thm}  \label{the:KKG_C_Ind}
  Let $G$, $H$ and~$D$ be as above.  Then there exist natural
  isomorphisms
  $$
  \KK^G_*(\C,\Ind_H^G D)\cong
  \K_{*+1}(\bvr(G,D_H)^H )\cong
  \K_{*+1}(\bvr(G,D)\cross H).
  $$
\end{thm}

\begin{proof}
  Since $\Comp(\ell^2 H)$ is finite dimensional, we have $\bvr(G,D_H)\cong
  \bvr(G,D)\otimes\Comp(\ell^2 H)$.  For any finite group~$H$ and any
  $H$\nbd{}$C^*$\brd{}algebra~$B$ there is a canonical isomorphism
  $\bigl(B\otimes\Comp(\ell^2 H)\bigr)^H \cong B \cross H$.
  Hence we have an isomorphism
  $$
  \bvr(G,D_H)^H \cong
  \bigl(\bvr(G,D)\otimes\Comp(\ell^2 H)\bigr)^H \cong
  \bvr(G,D)\cross H,
  $$
  which yields the second isomorphism of the theorem.  Actually, in
  the proof we will exclusively work with fixed point algebras.
  
  We now describe the cycles for $\KK^G_*(\C,\Ind_H^G D)$ more
  concretely, first for $*=1$.  These are given by pairs $(\Hilm,F)$
  where~$\Hilm$ is a $G$\nbd{}equivariant Hilbert module over
  $\Ind_H^G D$ and $F\in\Bound(\Hilm)$ satisfies
  \begin{equation}  \label{eq:KKG_C_cycle}
    F=F^*,\quad
    1-F^2\in\Comp(\Hilm),\quad
    gF-F\in\Comp(\Hilm)
  \end{equation}
  for all $g\in G$.  Since $\Ind_H^G D$ is a proper
  $G$\nbd{}$C^*$\brd{}algebra, the Equivariant Stabilization Theorem
  of~\cite{Meyer:KKG} applies.  Hence we can restrict attention to the
  case where~$\Hilm$ is the standard Hilbert module $\Ind_H^G
  D\otimes\ell^2(G\times\N)$ over $\Ind_H^G D$.  There is a natural
  isomorphism $\Ind_H^G D\otimes\ell^2(G\times\N) \cong \Ind_H^G
  \bigl(D\otimes\ell^2(H\times\N)\bigr)$.  Thus we may equivalently
  parametrize cycles for $\KK^G(\C,\Ind_H^G D)$ by operators on
  $\Ind_H^G (D\otimes\ell^2(H\times\N))$ that
  satisfy~\eqref{eq:KKG_C_cycle}.
  
  Elements of the Hilbert module $\Ind_H^G(D\otimes\ell^2(H\times\N))$
  are functions in $C_0(G,D\otimes\ell^2(H\times\N))$ that satisfy
  $f(gh)=\alpha_h \bigl( f(g)\bigr)$ for all $g\in G$, $h\in H$.  The
  right $\Ind_H^G D$\brd{}Hilbert module structure is given by
  pointwise multiplication and pointwise inner products.  The
  group~$G$ acts by left translation.  The operator~$F$ becomes a
  family of self-adjoint operators
  $$
  F =
  (F_g)_{g\in G},
  \qquad
  F_g\in\Bound(D\otimes\ell^2(H\times\N)) \cong
  \Mult\bigl(D_H\otimes\Comp\bigr).
  $$
  Since~$F$ preserves the covariance condition $f(gh)=\alpha_h \bigl(
  f(g)\bigr)$, we have
  \begin{equation}  \label{eq:F_covariant}
    F_{gh}= \alpha_h(F_g)
  \end{equation}
  for all $g\in G$, $h\in H$.  Secondly,
  $1-F_g^2\in\Comp(D\otimes\ell^2(H\times\N)) = D_H\otimes\Comp$, and
  the function $g\mapsto \norm{1-F_g^2}$ belongs to $C_0(G)$, whence,
  \begin{equation}  \label{eq:F_Fredholm}
    1 - F^2 \in C_0(G,D_H\otimes\Comp)
  \end{equation}
  Finally, $F_g-F_{xg}\in\Comp(D\otimes\ell^2(H\times\N))=
  D_H\otimes\Comp$ for all $g,x\in G$, and the function $g\mapsto
  \norm{F_g-F_{xg}}$ belongs to $C_0(G)$ for all $x\in G$.  Thus
  \begin{equation}  \label{eq:F_almostequivariant}
     F_g - F_{g'} \in\Comp(D\otimes\ell^2(H\times\N)) =
     D_H\otimes\Comp
  \end{equation}
  for all $g,g'\in G$, and
 \begin{equation}  \label{eq:F_vv}
    \lim_{g\to\infty} \norm{F_g-F_{xg}} = 0.
  \end{equation}
  for all $x\in G$.
  
  We claim that~\eqref{eq:F_vv} holds if and only if the function
  $g\mapsto F_g$ has vanishing variation in the sense of
  Definition~\ref{vanishingvariation}.  It is clear that vanishing
  variation implies~\eqref{eq:F_vv}.  Conversely, if we do not have
  vanishing variation, then there exist sequences $g_n, g_n'\to\infty$
  in~$G$ and $\epsilon>0$ such that the sequence $\{(g_n,g_n')\}$ is
  contained in some fixed entourage, but $\norm{F_{g_n} -
  F_{g_n'}}\ge\epsilon$ for all~$n$.  To say that $\{(g_n,g_n')\}$
  lies in an entourage means that $g_ng_n'^{-1}\in\Sigma$ for all~$n$,
  for some finite set $\Sigma\subset G$.  After extracting a
  subsequence, we may replace~$\Sigma$ by a singleton $\{x\}$, so that
  $g_n'=xg_n$ for all~$n$.  Hence~\eqref{eq:F_vv} is violated.  This
  proves that~\eqref{eq:F_vv} is equivalent to vanishing variation
  of~$F$.

  Equation~\eqref{eq:F_almostequivariant} means that~$F$ belongs to
  $\vvr(G,D_H)$.  Equation~\eqref{eq:F_covariant} means that
  $F\in\vvr(G,D_H)^H$.  In addition, we have $F=F^*$
  and~\eqref{eq:F_Fredholm}.  Let~$[F]$ be the image of~$F$ in
  $\bvr(G,D_H)$.  Equation~\eqref{eq:F_Fredholm} means that $P_F\defeq
  ([F]-1)/2$ is a projection in $\bvr(G,D_H)^H$.  Conversely, any
  projection in $\bvr(G,D_H)^H$ is of the form~$P_F$ for some
  cycle~$F$ for $\KK^G(\C,\Ind_H^G D)$.
  
  Since $\bvr(G,D_H)^H$ is matrix stable, we do not have to adjoin
  matrices to compute its $\K$\nbd{}theory.  The cycle~$F$ is
  degenerate if and only if $F_g=F_{xg}$ for all $g,x\in G$, that is,
  if and only if~$P_F$ is a \emph{constant} function on $G$.  The
  subalgebra of constant functions in $\vvr(G,D_H)^H$ is isomorphic to
  the stable multiplier algebra
  $\Mult\bigl(D_H\otimes\Comp\bigr)^H\cong\Mult(D\cross
  H\otimes\Comp)$ and hence has vanishing $\K$\nbd{}theory.  Two
  cycles $F$ and~$F'$ satisfy $P_F=P_{F'}$ if and only if~$F'$ is a
  compact perturbation of~$F$.  An operator homotopy between two
  cycles is the same as a homotopy between the associated projections.
  As a result, the equivalence relation generated by addition of
  degenerate cycles and operator homotopy for~$F$ is equivalent to the
  equivalence relation of stable homotopy equivalence of projections
  in $\bvr(G,D_H)^H$.  Since operator homotopy and homotopy generate
  the same equivalence relation on $\KK^G(\C,\Ind_H^G D)$, we obtain
  $\KK^G_1(\C,\Ind_H^G D)\cong\K_0(\bvr(G,D_H) ^H)$ as claimed.
  
  Consider now cycles for $\KK^G_0(\C,\Ind_H^G D)$.  Thus we also have
  a grading on our Hilbert module~$\Hilm$, and~$F$ is odd.  Since~$F$
  is self-adjoint and odd, knowing~$F$ is equivalent to knowing its
  restriction $U\colon \Hilm_+\to\Hilm_-$.  We may assume that the
  even and odd parts $\Hilm_\pm$ are isomorphic to
  $\Ind_H^G(D\otimes\ell^2(H\times\N))$.  As above, we obtain $U\in
  \vvr(G,D_H)^H $ with unitary image in $\bvr(G,D_H)^H$.  Conversely,
  any unitary in $\bvr(G,D_H)^H$ arises in this fashion.  We have
  similar criteria for degenerate Kasparov cycles, compact
  perturbations, and operator homotopy.  This yields
  $\KK^G_0(\C,\Ind_H^G D)\cong\K_1(\bvr(G,D_H)^H)$ as above.
\end{proof}

\section{The coarse co\brd{}assembly map}
\label{sec:coarse_coassembly}

In this section we define coarse co\brd{}assembly maps which we will
eventually identify with the analytic co\brd{}assembly maps occuring
in Theorem~\ref{bigtheorem2}.  We first introduce a
$\K$\nbd{}theoretic analogue of the coarse $\K$\nbd{}homology of a
space~$X$.  We do this equivariantly with respect to a finite group
action, and with coefficients.

Let~$X$ be a coarse space.  By our convention, (in particular by the
bounded geometry assumption), $X$ is coarsely equivalent to a discrete
coarse space.  Moreover, if~$H$ is a finite group acting isometrically
on~$X$, then~$X$ can even be arranged to be $H$\nbd{}equivariantly
coarsely equivalent to a discrete $H$\nbd{}space.  Since our
constructions depend only on the $H$\nbd{}equivariant coarse
equivalence class of~$X$, we may assume~$X$ discrete to begin with.
Let~$D$ be an $H$\nbd{}$C^*$\brd{}algebra and recall that~$D_H$
denotes $D\otimes\Comp(\ell^2 H)$ with its canonical $H$\nbd{}action.

Let~$E$ be an $H$\nbd{}invariant entourage of~$X$, and let $P_E(X)$ be
the simplicial complex whose vertices are the points of~$X$ and whose
simplices are the finite subsets $F\subseteq X$ for which $F\times F
\subset E$.  We do not distinguish between this simplicial complex and
its geometric realization, so that $P_E(X)$ is a locally compact
$H$\nbd{}space.  We view elements of~$P_E(X)$ as probability measures
on~$X$ in the usual fashion, so that the support of an element
of~$P_E(X)$ is a subset of~$X$.  We equip~$P_E(X)$ with the coarse
structure generated by the entourages
$$
\{(\mu,\nu)\in P_n(X) \times P_n(X) \mid
\supp(\mu)\times\supp(\nu)\subset F\},
$$
where~$F$ ranges over the entourages of~$X$.  It is clear that the
finite group~$H$ preserves this coarse structure.  If~$E$ is the
diagonal, we get back the space~$X$ itself.  If $E_1\subset E_2$, then
there is an injective, proper, $H$\nbd{}equivariant, continuous coarse
equivalence $P_{E_1}(X) \to P_{E_2}(X)$.  In particular, we obtain
$H$\nbd{}equivariant coarse equivalences $X\to P_E(X)$ for any~$E$
that contains the diagonal, and these maps are compatible with the
maps $P_{E_1}(X)\to P_{E_2}(X)$.

Since~$X$ is assumed separable, we can choose an increasing
sequence~$(E_n)$ of $H$\nbd{}invariant entourages of~$X$ such that any
entourage of~$X$ is contained in~$E_n$ for some~$n$.  Let $P_n(X)$
denote $P_{E_n}(X)$ and let $i_n\colon P_n(X) \to P_{n+1}(X)$ be the
canonical map.  These maps are injective and $H$\nbd{}equivariant.  We
denote the resulting inductive system by $\Rips(X)$.  Let
$\abs{\Rips(X)}$ be the inductive limit of this system, equipped with
the canonical topology.  We let
$$
\hat{C}_0(\Rips(X),D)\cross H \defeq
\varprojlim C_0(P_n(X),D)\cross H.
$$
This is a $\sigma$\nbd{}$C^*$\brd{}algebra in the sense
of~\cite{Phillips}.  Since the maps~$i_n$ are injective, the induced
maps on $C_0$\nbd{}functions are surjective.  Hence the maps
$\hat{C}_0(\Rips(X),D)\to C_0(P_n(X),D)$ are surjective for all
$n\in\N$.  This yields an isomorphism
$$
\hat{C}_0(\Rips(X),D) \cong
\{f\colon \abs{\Rips(X)}\to\C \mid
  \text{$f_{|_{P_n(X)}} \in C_0(P_{n}(X),D)$ for all $n\in\N$}
\}.
$$
Different choices of the sequence of entourages $(E_n)$ yield
isomorphic inductive systems $\Rips(X)$ and hence
$H$\nbd{}equivariantly isomorphic $\sigma$\nbd{}$C^*$\brd{}algebras
$\hat{C}_0(\Rips(X),D)$.  Even more, coarsely equivalent coarse spaces
yield homotopy equivalent systems of spaces $\Rips(X)$ and hence
homotopy equivalent $\sigma$\nbd{}$C^*$\brd{}algebras
$\hat{C}_0(\Rips(X),D)$ (see~\cite{EmersonMeyer}).  Hence the
following definition is legitimate.

\begin{defn}  \label{definitionofcoarsektheory}
  Let~$X$ be a coarse space with action of a finite group~$H$.  The
  \emph{$H$\nbd{}equivariant coarse $\K$\nbd{}theory of~$X$ with
  coefficents in~$D$} is defined by
  $$
  \KX^*_H(X,D) = \K_*\bigl(\hat{C}_0(\Rips(X),D)\cross H\bigr).
  $$
\end{defn}

The following lemma is proved in~\cite{EmersonMeyer}.

\begin{lemma}
  The assignment $X\mapsto \KX^*_H(X,D)$ is functorial from the
  category of coarse spaces with isometric actions of~$H$, and
  $H$\nbd{}equivariant coarse maps, to the category of Abelian groups
  and Abelian group homomorphisms.  That is, an equivariant coarse map
  $\phi\colon X\to Y$ induces a canonical homomorphism $\KX^*_H
  (Y,D)\to\KX^*_H(X,D)$ for any coefficient algebra $D$, and close
  maps induce the same map $\KX^*_H(Y,D)\to\KX^*_H(X,D)$.
\end{lemma}

We can now define the coarse co\brd{}assembly map.  Let~$X$ be a
coarse space equipped with an isometric action of a finite group $H$,
and let $D$ be an $H$\nbd{}$C^*$\brd{}algebra.  Construct the
inductive system $P_n(X)$ with injective, $H$\nbd{}equivariant coarse
equivalences $i_n\colon P_n(X)\to P_{n+1}(X)$ as above.  Since the
spaces $P_n(X)$ are themselves coarse spaces, we can form the algebras
$\vvr(P_n(X),D)\cross H$ and $\bvr(P_n(X),D)\cross H$.  By
functoriality of the various constructions involved, the maps~$i_n$
give rise to commutative diagrams with exact rows
$$
\xymatrix@C=1.2em{
  {C_0(P_{n+1}(X),D)\cross H} \ar[r] \ar@{->>}[d] &
  {\vvr(P_{n+1}(X),D)\cross H} \ar[r] \ar[d] &
  {\bvr(P_{n+1}(X),D)\cross H}  \ar[d]^{\cong} &
  \\
  {C_0(P_n(X),D)\cross H} \ar[r] &
  {\vvr(P_n(X),D)\cross H} \ar[r] &
  {\bvr(P_n(X),D)\cross H.} &
}
$$
We have surjective maps on kernel and quotient.  Hence the maps on
$\vvr(\cdots)$ are also surjective by the Snake Lemma.

\begin{lemma}[\cite{Phillips}]  \label{philipslemma}
  Suppose that $\alpha_n\colon A_{n+1} \to A_n$ is an inverse system
  of $C^*$\brd{}algebras with surjective maps~$\alpha_n$ for all~$n$.
  Let~$J_n$ be ideals in~$A_n$ such that the restriction of $\alpha_n$
  to $J_{n+1}$ maps $J_{n+1}$ surjectively onto~$J_n$.  Then the
  sequence
  $$
  0 \to
  \varprojlim J_n \to
  \varprojlim A_n \to
  \varprojlim A_n/J_n \to
  0
  $$
  is an exact sequence of $\sigma$\nbd{}$C^*$\brd{}algebras.
\end{lemma}

We have already introduced the inverse limit
$\hat{C}_0(\Rips(X),D)\cross H$ of the ideals $C_0(P_n(X),D)$ above.
Since the maps on $\bvr(P_n(X),D)\cross H$, are all isomorphisms, the
inverse limit of this system is again isomorphic to $\bvr(X,D)$.  Let
$$
\vvr(\Rips(X),D)\cross H\defeq
\varprojlim \vvr(P_n(X),D)\cross H.
$$
Phillips shows in~\cite{Phillips} how to associate to an exact
sequence of $\sigma$\nbd{}$C^*$\brd{}algebras a long exact sequence in
$\K$\nbd{}theory.  This is how we define the coarse co\brd{}assembly
map:

\begin{defn}  \label{def:coarse_coassembly}
  Let~$X$ be a coarse space, let~$H$ be a finite group acting on~$X$,
  and let~$D$ be an $H$\nbd{}$C^*$\brd{}algebra.  The
  \emph{$H$\nbd{}equivariant coarse co\brd{}assembly map with
  coefficients in~$D$} is the connecting map
  $$
  \mu^*_{D,H}\colon
  \K_{*+1}(\bvr(X,D)\cross H)\to\KX^*_H(X,D)
  $$
  associated to the exact sequence of
  $\sigma$\nbd{}$C^*$\brd{}algebras
  $$
  0 \to
  \hat{C}_0(\Rips(X),D)\cross H \to
  \vvr(\Rips(X),D)\cross H\to
  \bvr(X,D)\cross H \to
  0.
  $$
\end{defn}

Now let~$G$ be a countable discrete group, equipped with its canonical
coarse structure.  Let $(\Sigma_n)$ be an increasing sequence of
finite subsets of~$G$ with $G=\bigcup \Sigma_n$ and
$\Sigma_n=\Sigma_n^{-1}$ for all $n\in\N$ and $\Sigma_0=\{1\}$.  For
each~$n$, we get an entourage
$$
E_n \defeq \{(g_1,g_2) \in G\times G\mid g_1g_2^{-1} \in \Sigma_n\}.
$$
These entourages can be used to define the inductive system $\Rips(X)$
that occurs in the definition of the coarse $\K$\nbd{}theory of~$G$.
The entourage~$E_n$ is clearly invariant under right translation
by~$G$, so that the associated simplicial complex $P_n(G)= P_{E_n}(G)$
is a $G$\nbd{}space.  In addition, the action of~$G$ on $P_n(G)$ is
proper and $G$\nbd{}compact.  Hence $\abs{\Rips(G)}$ is a proper
$G$\nbd{}space.  It is not locally compact because the corresponding
simplicial complex is not locally finite.  Nevertheless, it is a model
for the universal proper $G$\nbd{}space~$\EG$ in the following sense:

\begin{prop}  \label{pro:Rips_EG}
  For any proper $G$\nbd{}space~$X$, the space of $G$\nbd{}equivariant
  maps $X\to\abs{\Rips(G)}$ is naturally homeomorphic to the space of
  cut-off functions on~$X$ and hence contractible.  Thus $\Rips(G)$ is
  a model for~$\EG$.
\end{prop}

\begin{proof}
  We view elements of $\abs{\Rips(G)}$ as functions $f\colon
  G\to[0,1]$ with compact support and $\sum f(g)=1$.  An equivariant
  map $h\colon X\to\abs{\Rips(G)}$ is already determined by the map
  $h_*\colon X\to[0,1]$, $h_*(x) = h(x)(1)$.  This map is continuous,
  satisfies $\sum h_*(gx)=1$ for all $x\in X$, and the support
  of~$h_*$ intersects each $G$\nbd{}compact subset of~$X$ in a compact
  subset.  That is, $h_*$ is a cut-off function on~$X$.  Conversely,
  any cut-off function arises from a unique map $X\to\abs{\Rips(G)}$.
  The space of cut-off functions on~$X$ is non-empty for any
  sufficiently regular, proper $G$\nbd{}space~$X$.  It is convex and
  hence contractible.  Hence there exists a unique map up to homotopy
  $X\to\abs{\Rips(G)}$ for any sufficiently regular, proper
  $G$\nbd{}space.  That is, $\abs{\Rips(G)}$ is universal.
\end{proof}

The above model for~$\EG$ is a variant of the model constructed by G.\
Kasparov and G.\ Skandalis in~\cite{KasparovSkandalis2}.  Let now~$X$
be any second countable, proper $G$\nbd{}space.  We let $(X_n)$ be an
increasing sequence of $G$\nbd{}compact subsets of~$X$ such that any
$G$\nbd{}compact subset of~$X$ is contained in~$X_n$ for some
$n\in\N$, and we let $\hat{C}_0(X)\defeq \hat{C}_0\bigl((X_n)\bigr)$.
The same reasoning as for $\hat{C}_0(\Rips(X))$ identifies
$$
\hat{C}_0(X) \cong
\{f\colon X\to\C \mid
  \text{$f_{|_Y}\in C_0(Y)$ for all $G$\nbd{}compact subsets
  $Y\subseteq X$}
\}.
$$
Clearly, any continuous, $G$\nbd{}equivariant map $X\to Y$ induces a
$*$\nbd{}homomorphism $\hat{C}_0(Y)\to\hat{C}_0(X)$.  Thus
equivariantly homotopic maps induce equivariantly homotopic
$*$\nbd{}homomorphisms.  Since the universal proper
$G$\nbd{}space~$\EG$ is uniquely determined up to $G$\nbd{}equivariant
homotopy equivalence, we obtain:

\begin{lemma}  \label{lem:Rips_EG}
  Let $H\subseteq G$ be a finite subgroup and let~$D$ be an
  $H$\nbd{}$C^*$\brd{}algebra.  Let~$\EG$ be any second countable,
  universal proper $G$\nbd{}space.  Then the
  $\sigma$\nbd{}$C^*$\brd{}algebras $\hat{C}_0(\EG,D)\cross H$ and
  $\hat{C}_0(\Rips(X),D)\cross H$ are homotopy equivalent.
\end{lemma}

\section{The Main Theorem}
\label{sec:main_theorem}

We are going to identify the coarse co\brd{}assembly map and the
analytic co\brd{}assembly map with appropriate coefficients.  We need
two preparatory results.  The first is well-known for
$C^*$\nbd{}algebras, and the proof for
$\sigma$\nbd{}$C^*$\brd{}algebras is exactly the same.

\begin{lemma}  \label{tinylemma}
  Let~$G$ be a discrete group and~$H$ a finite subgroup.  Let~$A$ be a
  $G$\nbd{}$\sigma$\nbd{}$C^*$\brd{}algebra and~$B$ an
  $H$\nbd{}$\sigma$\nbd{}$C^*$\brd{}algebra.
  \begin{enumerate}[(1)]
    \item The $\sigma$\nbd{}$C^*$\brd{}algebras $A\otimes\Ind_H^G B$
      and $\Ind_H^G(A\otimes B)$ are $G$\nbd{}equivariantly isomorphic.

   \item The $\sigma$\nbd{}$C^*$\brd{}algebras $\bigl(\Ind_H^G
     A\bigr)\cross G$ and $A\cross H$ are strongly Morita equivalent.

  \end{enumerate}
\end{lemma}

\begin{lemma}  \label{lem:induced_EG}
  Let~$G$ be a discrete group, $H$ a finite subgroup, and~$D$ an
  $H$\nbd{}$C^*$\brd{}algebra.  Then there is a canonical isomorphism
  $$
  \K_*\bigl(\hat{C}_0(\EG,D)\cross H\bigr) \cong
  \K_*\bigl(\hat{C}_0(\EG,\Ind_H^G D)\cross G \bigr).
  $$
\end{lemma}

\begin{proof}
  Lemma~\ref{tinylemma} implies
  \begin{multline*}
    \hat{C}_0(\EG,\Ind_H^G D) \cross G
    \cong
    (\hat{C}_0(\EG)\otimes\Ind_H^G D) \cross G
    \\ \cong
    \Ind_H^G( \hat{C}_0(\EG)\otimes D) \cross G
    \sim \hat{C}_0(\EG,D)\cross H,
  \end{multline*}
  where~$\cong$ denotes isomorphism and~$\sim$ denotes strong Morita
  equivalence.  The result follows.
\end{proof} 

\begin{lemma}  \label{lem:RKKG_as_Ktheory}
  Let~$G$ be a discrete group and~$X$ a second countable, proper
  $G$\nbd{}space, and let~$B$ be a $G$\nbd{}$C^*$\brd{}algebra.  Then
  there is a natural isomorphism
  $$
  \RKK^G_*(X;\C,B) \cong \K_*(\hat{C}_0(X,B)\cross G).
  $$
\end{lemma}

\begin{proof}
  We check that both groups agree on the level of cycles after some
  standard simplifications for cycles for $\RKK^G(X;\C,B)$.  Since
  $C_0(X,B)$ is a proper $G$\nbd{}$C^*$\brd{}algebra, the reduced and
  full crossed products for $C_0(X,B)$ are equal.  Moreover, the
  $C^*$\nbd{}categories of $G$\nbd{}equivariant Hilbert modules over
  $C_0(X,B)$ and of Hilbert modules over $C_0(X, B) \cross G$ are
  equivalent (see~\cite{Meyer:Fixed}).  That is, any
  $G$\nbd{}equivariant Hilbert module~$\Hilm$ over $C_0(X,B)$
  corresponds to a Hilbert module~$\tilde\Hilm$ over $C_0(X,B) \cross
  G$.  The correspondence is such that $\Bound(\tilde\Hilm)$ is
  naturally isomorphic to the $C^*$\nbd{}algebra $\Bound(\Hilm)^G$ of
  $G$\nbd{}equivariant, adjointable operators on~$\Hilm$.  The compact
  operators on~$\tilde\Hilm$ correspond to the \emph{generalized fixed
  point algebra} of $\Comp(\Hilm)$, which is the closed linear span of
  operators of the form $\sum_{g\in G} \alpha_g(\ket{\xi}\bra{\eta})$,
  where $\xi,\eta\in\Hilm$ are \emph{compactly supported} sections.
  (The support of~$\xi$ is the set of all $x\in X$ with $\xi_{x} \not=
  0$.)  More generally, if $T\in\Comp(\Hilm)$ has compact support in
  an appropriate sense, then $\sum_{g \in G} \alpha_g(T)$ belongs to
  the generalized fixed point algebra.
  
  The cycles for $\RKK^G(X;\C,B)$ are pairs $(\Hilm,F)$ where~$\Hilm$
  is a (graded) $G$\nbd{}equivariant Hilbert module over $C_0(X,B)$
  and where $F\in\Bound(\Hilm)$ satisfies $F=F^*$, $-1\le F\le 1$, $F$
  is odd in the graded case, $\phi(1-F^2)\in\Comp(\Hilm)$ for
  all $\phi\in C_0(X)$, and~$F$ is $G$\nbd{}equivariant.  We can
  arrange~$F$ to be exactly equivariant (see~\cites{Kasparov,
  Meyer:KKG}) since~$X$ is a proper $G$\nbd{}space.  By our category
  equivalence, all this data may be rewritten in terms of a pair
  $(\tilde\Hilm,\tilde{F})$, where~$\tilde\Hilm$ is a Hilbert module
  over $C_0(X,B)\cross G$ and $\tilde{F}\in\Bound(\tilde\Hilm)$
  satisfies $\tilde{F}=\tilde{F}^*$, $-1\le \tilde{F}\le 1$,
  and~$\tilde{F}$ is odd in the graded case.  There is an additional
  condition on $1-\tilde{F}^2$ which we now identify.

  For every $G$\nbd{}invariant, closed subset $Y\subseteq X$, we
  define restrictions of $\Hilm$ and~$\tilde\Hilm$ to~$Y$ by
  $$
  \Hilm_Y \defeq \Hilm/C_0(X\setminus Y)\cdot\Hilm,
  \qquad
  \tilde\Hilm_Y \defeq
  \tilde\Hilm/C_0\bigl((X\setminus Y)/G\bigr)\cdot\tilde\Hilm.
  $$
  We claim that $\phi(F^2-1)\in\Comp(\Hilm)$ for all $\phi\in C_0(X)$
  if and only if the operator on~$\tilde\Hilm_Y$ induced by
  $1-\tilde{F}^2$ is compact for all $G$\nbd{}compact subsets
  $Y\subseteq X$.

  Assume first that $\phi(F^2-1)\in\Comp(\Hilm)$ for all $\phi\in
  C_0(X)$.  Choose a $G$\nbd{}compact subset $Y\subseteq X$.  By the
  properness of the $G$\nbd{}action, there exists a function $\phi\in
  C_c(X)$ with $\sum_{g\in G} \phi(xg)=1$ for all $x\in Y$.  Then
  $1-F^2$ and $\sum_{g\in G} \alpha_g(\phi(1-F^2)) $ have the same
  restriction to~$Y$.  Since $\phi(1-F^2)$ is compact by hypothesis
  and has compact support, this sum belongs to $\Comp(\tilde\Hilm)$.
  Thus $(1-\tilde{F}^2)$ induces a compact operator on~$\tilde\Hilm_Y$
  for all $G$\nbd{}compact $Y\subseteq X$.

  Suppose conversely that $(1-\tilde{F}^2)$ induces a compact operator
  on~$\tilde\Hilm_Y$ for all $G$\nbd{}compact $Y\subseteq X$.  Let
  $\phi\in C_c(X)$.  Let~$Y$ be a $G$\nbd{}compact subset containing
  the support of~$\phi$.  Then $(1-\tilde{F}^2)|_Y$ belongs to the
  generalized fixed point algebra of~$\Hilm_Y$ and hence can be
  approximated by operators of the form $\sum_{g\in G} \alpha_g T$ for
  a finite rank operator~$T$ on~$\Hilm$ with compact support.  Hence
  the function $g\mapsto \phi\alpha_g T$ has compact support, so that
  $\phi \sum_{g\in G} \alpha_g T$ is a compact operator on~$\Hilm$.
  Since these operators approximate $\phi(1-F^2)$, the latter operator
  is also compact.

  Finally, we replace $C_0(X,B)$ by the
  $\sigma$\nbd{}$C^*$\brd{}algebra $\hat{C}_0(X,B)$.  A Hilbert
  module~$\tilde\Hilm$ over $C_0(X,B)\cross G$ automatically extends
  to a Hilbert module $\hat\Hilm\defeq \varprojlim \tilde\Hilm_{X_n}$
  over $\hat{C}_0(X,B)\cross G$.  Conversely, any Hilbert module over
  $\hat{C}_0(X,B)\cross G$ is of this form.  By definition, the
  algebra of compact operators on~$\hat\Hilm$ is $\varprojlim
  \Comp(\tilde\Hilm_{X_n})$.  Hence the above condition on
  $1-\tilde{F}^2$ is equivalent to $1-\tilde{F}^2\in\Comp(\hat\Hilm)$.
  We conclude that
  $$
  \RKK^G(X;\C,B)\cong
  \KK(\C,\hat{C}_0(X,B) \cross G) =
  \K(\hat{C}_0(X,B)\cross G),
  $$
  because these groups can be defined by exactly the same cycles and
  homotopies.
\end{proof}

\begin{thm}  \label{the:diagram}
  Let~$G$ be a discrete group, $H$ a finite subgroup of~$G$, and~$D$
  an $H$\nbd{}$C^*$\brd{}algebra.  Then there is a canonical
  isomorphism
  $$
  \KX^*_H(G,D) \cong \RKK^G(\EG;\C,\Ind_H^G D)
  $$
  such that the following diagram commutes:
  \begin{equation}
    \xymatrix{
      \KK^G_*(\C,\Ind_H^G D) \ar[d]^{\cong} \ar[r]^{p_\EG^*} &
      \RKK^G_*(\EG;\C,\Ind_H^G D) \ar[d]^{\cong} \\
      \K_{*+1}(\bvr(G,D)\cross H) \ar[r]^{\mu^*_{D,H}} &
      \KX^*_H(X,D).
    }
  \end{equation}
\end{thm}

\begin{proof}
  The lemmas \ref{lem:RKKG_as_Ktheory}, \ref{lem:induced_EG},
  and~\ref{lem:Rips_EG} yield isomorphisms
  \begin{multline*}
    \RKK^G(\EG;\C,\Ind_H^G D)
    \cong
    \K_*(\hat{C}_0(\EG,\Ind_H^G D)\cross G)
    \\ \cong
    \K_*(\hat{C}_0(\EG,D)\cross H)
    \cong
    \KX^*_H(G,D).
  \end{multline*}
  It remains to check that this makes our diagram commute.  We
  identify cycles for $\KK^G_1(\C,\Ind_H^G D)$ with certain elements
  $F\in\vvr(G,D)\cross H$ as in the proof of
  Theorem~\ref{the:KKG_C_Ind}.  Thus the isomorphism
  $\KK^G_1(\C,\Ind_H^G D) \to \K_*(\bvr(G,D)\cross H)$ maps~$[F]$ to
  the class of the projection $P_F\defeq([F]-1)/2$ in
  $\K_0(\bvr(G,D)\cross H)$.

  Let $0\to I\to E\to Q\to 0$ be an extension of
  $\sigma$\nbd{}$C^*$\brd{}algebras and let $P\in Q$ be a projection.
  Let $\partial\colon \K_0(Q)\to\K_1(I)$ be the connecting map.  To
  compute~$\partial[P]$, we lift~$P$ to any $\bar{P}\in E$ and embed
  $E\subseteq\Mult(I)$.  Then $\bar{P}^2-\bar{P}\in I$, and
  $\partial[P]\in \K_1(I)\cong\KK_1(\C,I)$ is represented by the cycle
  $(I,2\bar{P}-1)$.  We apply this to compute $\mu^*_{D,H}([P_F])$.
  It is easy to see that we can work with fixed point algebras instead
  of crossed products for the computation.  In particular, we have
  $\hat{C}_0(\EG,D_H\otimes\Comp)^H\cong
  \hat{C}_0(\EG,D\otimes\Comp)\cross H$.  Instead of lifting~$P_F$, we
  may just as well lift~$[F]$.  Thus, we have to extend
  $F\in\vvr(G,D_H)^H$ to~$\EG$.  Choose an $H$\nbd{}invariant
  continuous function $\phi\colon \EG\to\R_+$ such that (1) $S_Y\defeq
  \supp\phi \; \cap Y$ is compact for all $G$\nbd{}compact subsets
  $Y\subseteq\EG$, and (2) $\sum_{g\in G} \phi(xg)=1$ for all
  $x\in\EG$.  Let $L_Y\subseteq G$ be the set of all $g\in G$ with
  $S_Y\,g\cap S_Y\neq\emptyset$.  We let
  $$
  \bar{F}(x) \defeq \sum_{g\in G} \phi(xg)\,F(g^{-1}).
  $$
  If $x\in S_Y\,g$ for some $g\in G$, then $\bar{F}(x)$ is an average of
  $F(h)$ with $xh^{-1}\in S_Y$, so that $h\in L_Y^{-1}g$.  Hence
  $\norm{\bar{F}|_{S_Yg}- F_g}_\infty\to0$ for $g\to\infty$.  Thus
  $\bar{F}|_Y$ belongs to $\vvr(Y,D_H)^H $ for all $G$\nbd{}compact~$Y$, that
  is, $\bar{F}\in \vvr(\EG,D_H)^H$.  The quotient map $\varprojlim
  \vvr(\EG,D_H)^H\to\bvr(G,D_H)^H$ simply restricts a function on~$\EG$ to any
  $G$\nbd{}orbit in~$\EG$.  Hence $[\bar{F}]=[F]$ in $\bvr(G,D_H)^H$.
  Thus $\mu^*_{D,H}([P_F])$ is represented by the Kasparov
  cycle~$\bar{F}$ for $\KK(\C,\hat{C}_0(\EG,D_H\otimes\Comp)^H) \cong
  \KK(\C,\hat{C}_0(\EG,D)\cross H)$.
  
  Next we go around the diagram the other way.  By definition,
  $p_\EG^*(\alpha)$ is represented by the Kasparov cycle
  $(C_0(\EG,\Ind_H^G(D\otimes\ell^2(H\times\N)), F')$, where $F'
  f(x,g)\defeq F_g f(x,g)$ and~$G$ acts on $\EG\times G$ by $h\cdot
  (x,g)=(xh^{-1},hg)$.  The same formula makes sense if we replace
  $C_0$ by~$\hat{C}_0$.  To map this to $\K_0(\hat{C}_0(\EG,D)\cross
  H)$, we first have to make~$F'$ $G$\nbd{}equivariant.  We use the
  same cut-off function~$\phi$ as above to average~$F'$:
  $$
  F'' f(x,g)\defeq \sum_{h \in G} F_{h^{-1}g} \phi(xh)\, \cdot f(x,g)
  $$
  for all $x\in\EG$, $g\in G$.  This is a compact perturbation of~$F'$
  that is $G$\nbd{}equivariant.  Thus~$F''$ defines a multiplier of
  the generalized fixed point algebra of the Hilbert module
  $\hat{C}_0(\EG,\Ind_H^G(D\otimes\ell^2(H\times\N)))$.  Restriction
  to $\EG\times\{1\}\subseteq\EG\times G$ identifies this generalized
  fixed point algebra with $\hat{C}_0(\EG,D_H\otimes\Comp)^H\cong
  \hat{C}_0(\EG,D\otimes\Comp)\cross H$.  The composition of
  isomorphisms
  \begin{multline*}
    \RKK^G(\EG;\C,\Ind_H^G D)\cong
    \K_0(\hat{C}_0(\EG,\Ind_H^G D)\cross G)
    \\ \cong
    \K_0(\C,\hat{C}_0(\EG,D\otimes\Comp)\cross H) \defeq
    \KK_0(\C,\hat{C}_0(\EG,D\otimes\Comp)\cross H)
  \end{multline*}
  constructed above sends~$F'$ to the Kasparov cycle
  $(\hat{C}_0(\EG,D\otimes\Comp)\cross H, F''|_{\EG\times\{1\}})$.
  The reason is that Green's Imprimitivity Theorem is proved using the
  same manipulations of generalized fixed point algebras that we used
  above to view~$F''$ as a multiplier of
  $\hat{C}_0(\EG,D\otimes\Comp)\cross H$.  By construction,
  $F''|_{\EG\times\{1\}}=\bar{F}$.  This means that the diagram
  commutes on $\KK^G_1(\C,\Ind_H^G D)$.  The other parity can be
  treated by a similar argument, or by replacing~$D$ by $C_0(\R,D)$.
\end{proof}

Combining the above with Theorem~\ref{bigtheorem2} and
Proposition~\ref{pro:dual_Dirac_isomorphism}, we obtain:

\begin{thm}  \label{thebiggesttheoremofall1}
  Let~$G$ be a discrete, strongly geometrically finite group.  The
  following are equivalent:
  \begin{enumerate}[(1)]
    \item the $H$\nbd{}equivariant coarse co\brd{}assembly map with
      coefficients in~$D$ is an isomorphism for every finite
      subgroup~$H$ of~$G$ and every $H$\nbd{}$C^*$\brd{}algebra~$D$;

    \item the $H$\nbd{}equivariant coarse co\brd{}assembly map with
      coefficients in $C_0(\N)$ is an isomorphism for every finite
      subgroup~$H$ of~$G$.

    \item $G$ has a $\gamma$\nbd{}element.

  \end{enumerate}
  Moreover, (3) implies (1) and~(2) for an arbitrary discrete group.
\end{thm}

\begin{thm}  \label{thebiggesttheoremofall2}
  Let~$G$ be a discrete group with a $G$\nbd{}finite model for~$\EG$.
  Then~$G$ has a $\gamma$\nbd{}element if and only if the
  $H$\nbd{}equivariant coarse co\brd{}assembly map with trivial
  coefficients is an isomorphism for every finite subgroup~$H$ of~$G$.
\end{thm}

If~$G$ is torsion free, then the only finite subgroup is the trivial
group.  Since the coarse co\brd{}assembly map is an invariant of the
coarse structure, we obtain:

\begin{cor}
  If~$G$ is a torsion free, geometrically finite group, then the
  existence or non-existence of a $\gamma$\nbd{}element for~$G$ is a
  coarse invariant of~$G$.
\end{cor}

\end{document}